\documentclass[11pt,a4paper]{article}
\usepackage{amsmath}
\usepackage{amsfonts}

\newtheorem{theorem}{Theorem}

\newtheorem{corollary}[theorem]{Corollary}

\newtheorem{lemma}[theorem]{Lemma}

\newtheorem{proposition}[theorem]{Proposition}
\newtheorem{remark}[theorem]{Remark}

\newenvironment{proof}[1][Proof]{\noindent\textbf{#1.} }{\ \rule{0.5em}{0.5em}}
\input{tcilatex}
\begin{document}

\title{Boundary behaviour of Dirichlet series with applications to universal
series}
\date{}
\author{Stephen J. Gardiner and Myrto Manolaki}
\maketitle

\begin{abstract}
This paper establishes connections between the boundary behaviour of
functions representable as absolutely convergent Dirichlet series in a
half-plane\ and the convergence properties of partial sums of the Dirichlet
series on the boundary. This yields insights into the boundary behaviour of
Dirichlet series and Taylor series which have universal approximation
properties.
\end{abstract}

\section{General Dirichlet series\protect\footnotetext{%
\noindent 2010 \textit{Mathematics Subject Classification } 30B50, 30C85,
30K05, 30K10, 31A15.}}

We consider a general Dirichlet series of the form $f(s)=\sum\limits_{n=1}^{%
\infty }a_{n}e^{-\lambda _{n}s}$, where $(\lambda _{n})$ is an unbounded
strictly increasing sequence of nonnegative real numbers. If this series
converges when $s=s_{0}$, then, as is well known, it converges uniformly
throughout any angular region of the form $\{s\in \mathbb{C}:\left\vert \arg
(s-s_{0})\right\vert <\pi /2-\delta \}$, where $\delta \in (0,\pi /2)$. In
particular, $f$ is defined on the half-plane $\{\func{Re}s>\func{Re}s_{0}\}$
and has a nontangential limit at $s_{0}$. No such conclusion may be drawn if
we know merely that some subsequence of the partial sums of the series
converges at $s=s_{0}$ (see, for example, Bayart \cite{Bay}). Nevertheless,
we show below that there are strong connections between the nontangential
boundary behaviour of $f$ and the limiting behaviour of such subsequences.
This has implications for the boundary behaviour of universal Dirichlet
series and universal Taylor series.

A typical point of the complex plane will be written as $s=\sigma +it$. We
denote by $\mathbb{C}_{+}$ the right-hand half-plane $\{\sigma >0\}$, and by 
$D(\mathbb{C}_{+})$ the space of all holomorphic functions on $\mathbb{C}%
_{+} $ which can be represented there as an absolutely convergent Dirichlet
series $\sum\limits_{n=1}^{\infty }a_{n}e^{-\lambda _{n}s}$. Such a
representation is unique (see Theorem 6 of \cite{HR}). We define the partial
sum%
\begin{equation*}
S_{m}(s)=\sum_{n=1}^{m}a_{n}e^{-\lambda _{n}s}\text{ \ \ \ }(s\in \mathbb{C}%
),
\end{equation*}%
and denote by $\mathrm{nt}\lim_{s\rightarrow \zeta }f(s)$ the nontangential
limit, when it exists (finitely), of $f$ at a point $\zeta $ of the
imaginary axis $i\mathbb{R}$.

\begin{theorem}
\label{one}Let $f\in D(\mathbb{C}_{+})$, let $(S_{m_{k}})$ be a subsequence
of $(S_{m})$, and define%
\begin{equation*}
E=\{\zeta \in i\mathbb{R}:S(\zeta ):=\lim_{k\rightarrow \infty
}S_{m_{k}}(\zeta )\text{ exists}\}
\end{equation*}%
and%
\begin{equation*}
F=\{\zeta \in i\mathbb{R}:f(\zeta ):=\mathrm{nt}\lim_{s\rightarrow \zeta
}f(s)\text{ exists}\}.
\end{equation*}%
Then $f=S$ almost everywhere (Lebesgue) on $E\cap F$.
\end{theorem}

Theorem \ref{one} may be applied to the Riemann zeta function $\zeta $,
since the function $s\mapsto \zeta (s+1)$ belongs to $D(\mathbb{C}_{+})$. It
tells us that, if a subsequence $\left( \sum_{1}^{m_{k}}n^{-s}\right) $ of
the Dirichlet polynomials is pointwise convergent on a set $E\subset
\{\sigma =1\}$, then the limit function must equal $\zeta (s)$ almost
everywhere on $E$.

\bigskip

\noindent \textbf{Definition }By a \textit{fat approach region} (relative to 
$\mathbb{C}_{+}$) to a point $it_{0}$ of the imaginary axis, we mean a set
of the form 
\begin{equation*}
\Omega (it_{0})=\{\sigma +it:\left\vert t-t_{0}\right\vert <a\text{ \ and \ }%
\phi (t-t_{0})<\sigma <b\},
\end{equation*}%
where $a,b>0$ and $\phi :[-a,a]\rightarrow \lbrack 0,\infty )$ is a
Lipschitz function satisfying 
\begin{equation*}
\int_{-a}^{a}y^{-2}\phi (y)dy<\infty .
\end{equation*}

It is easy to see that such a set contains Stolz regions at $it_{0}$ of
arbitrarily large aperture, that is, $\mathrm{nt}\lim_{s\rightarrow
it_{0}}\chi _{\Omega (it_{0})}(s)=1$, where $\chi _{A}$ denotes the
characteristic function of a set $A$. Examples of fat approach regions to $0$
are the half-discs $\{s:\left\vert s-a\right\vert <a,\sigma <a\}$, where $%
a>0 $, and\ sets of the form $\{\sigma +it:\left\vert t\right\vert ^{\alpha
}<\sigma <1\}$, where $\alpha >1$.

\bigskip

\begin{theorem}
\label{two}Let $f\in D(\mathbb{C}_{+})$, let $h$ be a positive harmonic
function on $\mathbb{C}_{+}$, and suppose that $e^{-h}f$ is bounded on two 
\textit{fat approach region}s, $\Omega (it_{1})$ and $\Omega (it_{2})$,
where $t_{1}<t_{2}$. If a subsequence $(S_{m_{k}})$ is uniformly bounded on
an open interval of $i\mathbb{R}$ containing $it_{1}$ and $it_{2}$, then $f$
is bounded on any rectangle of the form $\{\sigma +it:0<\sigma
<1,t_{1}+\varepsilon <t<t_{2}-\varepsilon \}$, where $\varepsilon >0$.
\end{theorem}

\begin{remark}
\label{R0}(a)\ In Theorem \ref{two} the function $f$ must then have a
nontangential limit $f(\zeta )$ at almost every point $\zeta $ of the
interval $(it_{1},it_{2})$, by the local version of the Fatou theorem. Thus,
by Theorem \ref{one}, $f(\zeta )=\lim_{k\rightarrow \infty }S_{m_{k}}(\zeta
) $ at almost every point of $(it_{1},it_{2})$ where the latter limit exists.%
\newline
(b) In the special case where the function $f$ itself is bounded on two 
\textit{fat approach region}s, $\Omega (it_{1})$ and $\Omega (it_{2})$, the
proof of Theorem \ref{two} shows that $(S_{m_{k}})$ is uniformly bounded on
the rectangle $\{\sigma +it:0<\sigma <1,t_{1}<t<t_{2}\}$, and so $f$ is also
bounded there.
\end{remark}

The High Indices Theorem of Hardy and Littlewood (see \cite{HL} or Theorem
23.1 of \cite{Kor}) says that, if $\inf_{n}\lambda _{n+1}/\lambda _{n}>1$
and $f(it_{0}):=\lim_{\sigma \rightarrow 0+}f(\sigma +it_{0})$\ exists, then 
$S_{n}(it_{0})\rightarrow f(it_{0})$. We will now see a complementary result
concerning Dirichlet series which possess a different type of gap structure.

If $f$ is a holomorphic function on $\mathbb{C}_{+}$ which has bounded
derivative on a fat approach region to $\zeta \in i\mathbb{R}$, it is easy
to see that $f$ has a finite nontangential limit there, which we again
denote by $f(\zeta )$ in the next result.

\begin{theorem}
\label{three}Let $f\in D(\mathbb{C}_{+})$ and suppose that $f^{\prime }$ is
bounded on a fat approach region to $\zeta \in i\mathbb{R}$. If a
subsequence $(S_{m_{k}})$ is uniformly bounded on an open interval of $i%
\mathbb{R}$ containing $\zeta $, and $\lambda _{m_{k}+1}/\lambda
_{m_{k}}\rightarrow \infty $, then $S_{m_{k}}(\zeta )\rightarrow f(\zeta )$.
\end{theorem}

Theorem \ref{three} fails if we drop the requirement that $\lambda
_{m_{k}+1}/\lambda _{m_{k}}\rightarrow \infty $. To see this, let $%
f(s)=(e^{s}-1)^{-1}$. Then $f\in D(\mathbb{C}_{+})$, since $%
f(s)=\sum_{1}^{\infty }e^{-ns}$, and $f^{\prime }$ is bounded on $\{s\in 
\mathbb{C}_{+}:\left\vert s-i\pi \right\vert <\pi /2\}$. However, if we take 
$\zeta =i\pi $ and $\lambda _{m_{k}}=m_{k}=2k$, then%
\begin{equation*}
\left\vert S_{m_{k}}(it)\right\vert =\left\vert \frac{1-e^{-m_{k}it}}{%
e^{it}-1}\right\vert \leq \frac{2}{\sqrt{2}}=\sqrt{2}\text{ \ \ \ }(\pi
/2<t<3\pi /2)
\end{equation*}%
and $S_{m_{k}}(\zeta )=0$ for all $k\in \mathbb{N}$, yet $f(\zeta )=-1/2$.

Theorems \ref{one} and \ref{two} provide analogues for Dirichlet series of
results recently established for Taylor series in \cite{GM} and \cite{G14},
respectively. Theorem \ref{three} implies a corresponding result for Taylor
series $\sum a_{n}z^{n}$ in the unit disc $\mathbb{D}$ with (pure) Ostrowski
gaps, which can readily be deduced using the substitution $z=e^{-s}$ (see
Corollary \ref{C6} in Section \ref{pfapps}). In the next section we present
applications of the above theorems to universal Dirichlet series and
universal Taylor series.

\section{Applications to universal series}

We now focus on ordinary Dirichlet series, where $\lambda _{n}=\log n$. We
say that a function $f$ in $D(\mathbb{C}_{+})$ belongs to the collection $%
UD_{0}$, if it is of the form $f(s)=\sum\limits_{n=1}^{\infty }\dfrac{a_{n}}{%
n^{s}}$ and, for any continuous function $g:i\mathbb{R}\rightarrow \mathbb{C}
$ and any $b>0$, there is a subsequence $(S_{m_{k}})$ of the partial sums
that converges to $g$ uniformly on $[-bi,bi]$. Since $UD_{0}$ contains all
Dirichlet series with universal approximation properties in $\mathbb{C}%
\backslash \mathbb{C}_{+}$ as described by Bayart \cite{Bay}, it follows
from Theorem 6 of that paper that membership of $UD_{0}$ is topologically
generic in $D(\mathbb{C}_{+})$. (The topology is defined using the family of
seminorms $\left\Vert \sum a_{n}n^{-s}\right\Vert _{\sigma }=\sum \left\vert
a_{n}\right\vert n^{-\sigma },$ $\sigma >0$.)\ Further results on universal
Dirichlet series may be found in \cite{BGNP}, \cite{DM}, \cite{KNSS}, and 
\cite{NP}. Theorems \ref{one} and \ref{two} yield the following strong
information about the boundary behaviour of functions in $UD_{0}$, and so
about universal Dirichlet series in particular.

\begin{corollary}
\label{C1}Let $f\in UD_{0}$. Then, for almost every $\zeta \in i\mathbb{R}$,
the set $f(\Gamma )$ is dense in $\mathbb{C}$ for every open isosceles
triangle $\Gamma \subset \mathbb{C}_{+}$ with vertex at $\zeta $ and
horizontal axis of symmetry.
\end{corollary}

\begin{corollary}
\label{C2}Let $f\in UD_{0}$. Then, for any disc $D$ centred on $i\mathbb{R}$%
, the set $\mathbb{C}\backslash f(D\cap \mathbb{C}_{+})$ has zero
logarithmic capacity.
\end{corollary}

\begin{corollary}
\label{C3}Let $f\in UD_{0}$. Then there is a residual set $Z\subset \mathbb{R%
}$ such that $\{f(\sigma +it):0<\sigma <1\}$ is dense in $\mathbb{C}$ for
every $t\in Z$.
\end{corollary}

\begin{corollary}
\label{C4}Let $f\in D(\mathbb{C}_{+})$ and $h$ be a positive harmonic
function on $\mathbb{C}_{+}$, and suppose that $e^{-h}f$ is bounded on two 
\textit{fat approach region}s $\Omega (it_{1})$ and $\Omega (it_{2})$, where 
$t_{1}\neq t_{2}$. Then $f\not\in UD_{0}$.
\end{corollary}

The above corollaries are inspired by results recently established for
universal Taylor series in \cite{G13} and \cite{G14}. The hypotheses of
Corollary \ref{C4} clearly hold if, for example, $f$ is bounded on two discs
in $\mathbb{C}_{+}$ that are tangent to distinct points of $i\mathbb{R}$.

Next we use Theorem \ref{three} to give a new result for Taylor series. We
recall that a holomorphic function $f$ on $\mathbb{D}$ is said to possess a
universal Taylor series if, for any compact set $K\subset \mathbb{D}^{c}$
with connected complement and any continuous function $g:K\rightarrow 
\mathbb{C}$ that is holomorphic on $K^{\circ }$, there is a subsequence $%
(T_{m_{k}})$ of the partial sums of the Taylor series of $f$ about $0$ which
converges to $g$ uniformly on $K$. The collection of all such functions is
denoted by $U(\mathbb{D},0)$. Nestoridis \cite{N96} has shown that this
universal approximation property is topologically generic for holomorphic
functions on $\mathbb{D}$ (with respect to the topology of local uniform
convergence). Let $\omega \subset \mathbb{D}$ be open and $w\in \partial 
\mathbb{D}$. We call $\omega $ a fat approach region to $w$ (relative to $%
\mathbb{D}$) if the set $\{1-w^{-1}z:z\in \omega \}$ is a fat approach
region to $0$ (relative to $\mathbb{C}_{+}$).

\begin{corollary}
\label{C5}Let $f$ be a holomorphic function on $\mathbb{D}$. If $f^{\prime }$
is bounded on a fat approach region to some point of $\partial \mathbb{D}$,
then $f\not\in U(\mathbb{D},0)$.
\end{corollary}

Corollary \ref{C5} complements a result of Costakis and Katsoprinakis
(Theorem 2 of \cite{CK}) concerning the behaviour of antiderivatives of
universal Taylor series on internally tangential discs.

\section{Proof of Theorem \protect\ref{one}}

Our proof of Theorem \ref{one}\ relies on the following convergence result
for harmonic measures, taken from \cite{GM}. Although it is valid more
generally, we will only need the case where $\Omega $ is a plane domain
possessing a Green function $G_{\Omega }(\cdot ,\cdot )$. For any (nonempty)
open set $\omega $, any Borel set $A$ and any point $z$ in $\mathbb{\omega }$%
, we denote by $\mu _{z}^{\omega }(A)$ the harmonic measure of $A\cap
\partial \omega $ for $\omega $ evaluated at $z$. We refer to the book \cite%
{AG} for potential theoretic concepts.

\begin{theorem}
\label{tool}Let $\xi _{0}\in \Omega $ and $\omega $ be an open subset of $%
\Omega $. Suppose that $(v_{k})$ is a decreasing sequence of subharmonic
functions on $\omega $ such that $v_{1}/G_{\Omega }(\xi _{0},\cdot )$ is
bounded above and $\lim_{k\rightarrow \infty }v_{k}<0$ on $\omega $. If $\mu
_{z_{1}}^{\omega }(\partial \Omega )>0$ for some $z_{1}$, then $\mu
_{z_{1}}^{\{v_{k}<0\}}(\partial \Omega )>0$ for all sufficiently large $k$.
\end{theorem}

We will also repeatedly use the following lemmas. As usual, the positive and
negative parts of a real number $t$ will be denoted by $t^{+}$ and $t^{-}$,
respectively.

\begin{lemma}
\label{L0}Let $f(s)=\sum\limits_{n=1}^{\infty }a_{n}e^{-\lambda _{n}s}$,
where $f\in D(\mathbb{C}_{+})$, and let $u_{m}=\lambda _{m+1}^{-1}\log
|S_{m}-f|$. Then, for any $\varepsilon >0$, there exists $m_{0}\in \mathbb{N}
$ such that 
\begin{equation*}
u_{m}(s)\leq -\sigma +\varepsilon \text{ \ \ \ }(s=\sigma +it,\sigma
>\varepsilon ,m\geq m_{0}).
\end{equation*}
\end{lemma}

\begin{proof}
Given $\varepsilon >0$, we have%
\begin{eqnarray*}
u_{m}(s) &=&\frac{\log \left\vert \sum_{m+1}^{\infty }a_{n}e^{-\lambda
_{n}s}\right\vert }{\lambda _{m+1}} \\
&\leq &\frac{\log \left( e^{\lambda _{m+1}(\varepsilon -\sigma
)}\sum_{m+1}^{\infty }\left\vert a_{n}\right\vert e^{-\lambda
_{n}\varepsilon }\right) }{\lambda _{m+1}}\leq \varepsilon -\sigma \text{ \
\ \ \ \ }(\sigma >\varepsilon ),
\end{eqnarray*}%
for all sufficiently large $m$, since $f\in D(\mathbb{C}_{+})$.
\end{proof}

\begin{lemma}
\label{L}Let $f(s)=\sum\limits_{n=1}^{\infty }a_{n}e^{-\lambda _{n}s}$,
where $f\in D(\mathbb{C}_{+})$, let $\sigma _{0}<0$ and let $K\subset
\{\sigma _{0}<\sigma \leq 0\}$ be a nonpolar compact set. Then there exists $%
c>0$ such that%
\begin{equation*}
\left\vert S_{m}\right\vert \leq \max \left\{ \sup_{K}\left\vert
S_{m}\right\vert ,1\right\} e^{c\lambda _{m}G_{\mathbb{C}\backslash
K}(\sigma _{0}-1,\cdot )}\text{ \ on \ }\{\sigma >\sigma _{0}\}\backslash K%
\text{ \ for all }m\text{.}
\end{equation*}
\end{lemma}

\begin{proof}
We dismiss the trivial case where $m=1$ and $\lambda _{1}=0$. Since $f\in D(%
\mathbb{C}_{+})$, we know (cf. Theorem 8.1 of \cite{Ap}) that%
\begin{equation*}
\underset{m\rightarrow \infty }{\lim \sup }\frac{\log \left(
\sum_{1}^{m}\left\vert a_{n}\right\vert \right) }{\lambda _{m}}\leq 0.
\end{equation*}%
Also, 
\begin{equation*}
\left\vert S_{m}(s)\right\vert \leq \sum_{n=1}^{m}\left\vert
a_{n}\right\vert e^{-\lambda _{n}\sigma }\leq e^{\lambda _{m}\sigma
^{-}}\sum\limits_{n=1}^{m}\left\vert a_{n}\right\vert .
\end{equation*}%
Thus the functions 
\begin{equation*}
s\mapsto \left( \dfrac{\log \left\vert S_{m}(s)\right\vert }{\lambda _{m}}%
-\sigma ^{-}\right) ^{+}\text{ \ \ \ }(m\in \mathbb{N})
\end{equation*}%
are uniformly bounded on $\mathbb{C}$ and converge uniformly to $0$ as $%
m\rightarrow \infty $.

Let $b_{m}=\max \left\{ \sup_{K}\left\vert S_{m}\right\vert ,1\right\} $. By
the previous paragraph there exists a positive constant $c_{1}$ such that 
\begin{equation*}
\frac{\log (\left\vert S_{m}\right\vert /b_{m})}{\lambda _{m}}\leq \frac{%
\log (\left\vert S_{m}\right\vert )}{\lambda _{m}}\leq c_{1}\text{ \ \ \ \
on }\{\sigma +it:\sigma \geq \sigma _{0}\}\text{ \ for all }m\text{.}
\end{equation*}%
Also, $G_{\mathbb{C}\backslash K}(\sigma _{0}-1,\cdot )$ has a positive
limit at $\infty $. Thus, applying the maximum principle in $\{\sigma
>\sigma _{0}\}\backslash K$, we see that there is a positive constant $c$
such that 
\begin{equation*}
\log \left( \frac{\left\vert S_{m}\right\vert }{b_{m}}\right) \leq c\lambda
_{m}G_{\mathbb{C}\backslash K}(\sigma _{0}-1,\cdot )\text{ \ \ on }\{\sigma
>\sigma _{0}\}\backslash K,
\end{equation*}%
as required.
\end{proof}

\bigskip

\begin{proof}[Proof of Theorem \protect\ref{one}]
We adapt an argument from \cite{GM}. Let $f(s)=\sum\limits_{n=1}^{\infty
}a_{n}e^{-\lambda _{n}s}$, where $f\in D(\mathbb{C}_{+})$, and let $%
(S_{m_{k}})$, $E$ and $F$ be as in the statement of Theorem \ref{one}. For
each $\zeta \in i\mathbb{R}$, let $\Gamma (\zeta )$ denote the open
triangular region with vertices $\{\zeta ,\zeta +1\pm i\}$. We suppose, for
the sake of contradiction, that the conclusion of the theorem fails. Then,
first multiplying $f$ by a suitable unimodular constant, we may fix a number 
$a\geq 1$ large enough so that the set 
\begin{equation*}
\{\zeta \in E\cap F:\func{Re}\left( S(\zeta )-f(\zeta )\right) \geq a^{-1}%
\text{, }\left\vert f\right\vert \leq a\text{ on }\Gamma (\zeta ),\text{ }%
\left\vert S_{m_{k}}(\zeta )\right\vert \leq a\text{ for all }k\}
\end{equation*}%
has positive Lebesgue measure. Let $K$ be a compact subset of the above set,
also of positive measure, let $\Omega =\mathbb{C}\backslash K$ and $\omega
=\cup _{\zeta \in K}\Gamma (\zeta )$. (We can choose $K$ so that $\omega $
is connected.)

We define the subharmonic functions 
\begin{equation*}
u_{k}=\frac{1}{\lambda _{m_{k}+1}}\log \left( \frac{\left\vert
S_{m_{k}}-f\right\vert }{2a}\right) \text{ \ \ \ on \ }\mathbb{C}_{+}\text{
\ \ \ }(k\in \mathbb{N}).
\end{equation*}%
By Lemma \ref{L0} we may choose a sequence $\sigma _{k}\downarrow 0$ such
that%
\begin{equation}
u_{j}(s)\leq -\sigma /2\text{ \ \ \ }(\sigma _{k}\leq \sigma \leq 1,j\geq k).
\label{luc}
\end{equation}%
Also,%
\begin{equation}
u_{k}\leq \frac{1}{\lambda _{m_{k}}}\log ^{+}\left( \frac{\max \{\left\vert
S_{m_{k}}\right\vert ,\left\vert f\right\vert \}}{a}\right) \leq cG_{\Omega
}(-2,\cdot )\text{ \ on \ }\mathbb{\omega }\cup K\text{ \ }(k\geq m_{0}),
\label{uk}
\end{equation}%
by Lemma \ref{L} (where $\sigma _{0}=-1$, say). We now define $v_{k}$ to be
the solution, $H_{\psi _{k}}^{\omega }$, to the Dirichlet problem in $\omega 
$ with boundary data $\psi _{k}$, where%
\begin{equation*}
\psi _{k}(s)=\left\{ 
\begin{array}{cc}
-\sigma /2 & \text{if }\sigma \geq \sigma _{k} \\ 
cG_{\Omega }(-2,\cdot ) & \text{if }0<\sigma <\sigma _{k} \\ 
0 & \text{if }\mathbb{\sigma }=0%
\end{array}%
\right. .
\end{equation*}%
Then $u_{k}\leq v_{k}$ on $\omega $, by (\ref{luc}) and (\ref{uk}), and $%
(v_{k})$ is a decreasing sequence of harmonic functions on $\omega $ with
limit $-\sigma /2$ \ on $\partial \omega $.

Since $K$ has positive linear measure, we know that $\mu _{z}^{\omega }(K)>0$
when $z\in \omega $. (This follows from the F. and M. Riesz theorem and the
fact that $\partial \omega $ is rectifiable; see Section VI.1 of \cite{GaMa}%
.) Thus, by Theorem \ref{tool}, there exists $k^{\prime }\in \mathbb{N}$
such that the open set $\omega _{1}:=\omega \cap \{v_{k^{\prime }}<0\}$ is
nonempty and 
\begin{equation}
\mu _{w}^{\omega _{1}}(\partial \Omega )>0\text{ \ for some }w\in \omega
_{1}.  \label{ne}
\end{equation}%
When $k\geq k^{\prime }$ we have $u_{k}<0$ on $\omega _{1}$, so $\left\vert
S_{m_{k}}-f\right\vert <2a$, and hence $\left\vert S_{m_{k}}\right\vert <3a$%
, on $\omega _{1}$. Now $S_{m_{k}}=H_{S_{m_{k}}}^{\omega _{1}}$on $\omega
_{1}$, so by dominated convergence 
\begin{equation*}
f=H_{\phi }^{\omega _{1}}\text{ \ on \ }\omega _{1}\text{, \ \ where \ \ }%
\phi =\left\{ 
\begin{array}{cc}
f & \text{ \ on \ }\partial \omega _{1}\cap \mathbb{C}_{+} \\ 
S & \text{ \ on \ }\partial \omega _{1}\cap i\mathbb{R}\subset K%
\end{array}%
\right. .
\end{equation*}%
However, we also know that 
\begin{equation*}
f=H_{f}^{\omega }=H_{H_{f}^{\omega }}^{\omega _{1}}=H_{f}^{\omega _{1}}\text{
\ on \ }\omega _{1}
\end{equation*}%
(see Theorem 6.3.6\ of \cite{AG}). Thus $H_{(S-f)\chi _{K}}^{\omega
_{1}}\equiv 0$. This yields a contradiction, since our choice of $K$ and $%
\omega _{1}$ ensures that%
\begin{equation*}
\func{Re}H_{(S-f)\chi _{K}}^{\omega _{1}}(w)\geq a^{-1}\mu _{w}^{\omega
_{1}}(K)>0.
\end{equation*}%
Theorem \ref{one} is now established.
\end{proof}

\section{Proof of Theorem \protect\ref{two}}

We will need the following variant of Proposition 4 in \cite{GM}.

\begin{proposition}
\label{P}Let $\omega $ be a fat approach region to $it_{0}\in i\mathbb{R}$
relative to $\mathbb{C}_{+}$. Further, let $(u_{k})$ be a sequence of
subharmonic functions on $\omega $ such that $u_{k}(s)\leq b\sigma $ for all 
$k$, where $b\geq 0$, and $\lim \sup_{k\rightarrow \infty }u_{k}<0$. Then
there exists $k_{0}\in \mathbb{N}$ such that%
\begin{equation*}
\mathrm{nt~}\underset{s\rightarrow it_{0}}{\lim \sup }\frac{\sup_{k\geq
k_{0}}u_{k}(s)}{\sigma }<0.
\end{equation*}
\end{proposition}

\begin{proof}[Proof of Proposition]
We choose another fat approach region $\omega _{1}$ to $it_{0}$ such that $%
\overline{\omega }_{1}\backslash \{it_{0}\}\subset \omega $, and observe (by
Corollary 5.7.2 of \cite{AG}) that there are strictly decreasing sequences $%
\sigma _{k}\rightarrow 0$ and $\tau _{k}\rightarrow 0$\ such that 
\begin{equation*}
u_{j}(s)\leq -\tau _{k}\text{ \ \ \ }(s\in \partial \omega _{1},\sigma \geq
\sigma _{k},j\geq k).
\end{equation*}%
Let $v_{k}=H_{\psi _{k}}^{\omega _{1}}$, where%
\begin{equation*}
\psi _{k}(s)=\left\{ 
\begin{array}{cc}
-\tau _{1} & (\sigma \geq \sigma _{1}) \\ 
-\tau _{j} & (\sigma _{j}\leq \sigma <\sigma _{j-1};j=2,...,k) \\ 
b\sigma & (0\leq \sigma <\sigma _{k})%
\end{array}%
\right. .
\end{equation*}%
Then $u_{j}\leq v_{k}$ on $\omega _{1}$ whenever $j\geq k$, and $(v_{k})$ is
a decreasing sequence of harmonic functions on $\omega _{1}$ with $\lim
v_{k}<0$ there.

We now recall (see Chapter 9 of \cite{AG}) that a set $E\subset \mathbb{C}%
_{+}$ is said to be minimally thin at $it_{0}$ (with respect to $\mathbb{C}%
_{+}$) if there is a positive superharmonic function $h$ on $\mathbb{C}_{+}$
such that%
\begin{equation*}
\inf_{E}\frac{h}{P_{it_{0}}}>\inf_{\mathbb{C}_{+}}\frac{h}{P_{it_{0}}},
\end{equation*}%
where $P_{it_{0}}$ denotes the Poisson kernel for $\mathbb{C}_{+}$ with pole
at $it_{0}$; that is, 
\begin{equation*}
P_{it_{0}}(s)=\frac{\sigma }{\sigma ^{2}+(t-t_{0})^{2}}\text{ \ \ \ }%
(s=\sigma +it\in \mathbb{C}_{+}).
\end{equation*}%
Theorem 9.7.1 of \cite{AG} \ tells us that, because $\omega _{1}$ is a fat
approach region at $it_{0}$, the set $\mathbb{C}_{+}\backslash \omega _{1}$
is minimally thin at $it_{0}$. We can now apply Proposition 4 of \cite{GM}
to see that there exist $k_{0}\in \mathbb{N}$ and a set $E\subset \mathbb{C}%
_{+}$, minimally thin at $it_{0}$, such that 
\begin{equation*}
\underset{s\rightarrow it_{0},s\in \mathbb{C}_{+}\backslash E}{\lim }\frac{%
v_{k_{0}}(s)}{\sigma }<0.
\end{equation*}%
Next, we appeal to Theorem 6 of \cite{BD} to see that%
\begin{equation*}
\mathrm{nt~}\underset{s\rightarrow it_{0}}{\lim }\frac{v_{k_{0}}(s)}{\sigma }%
<0.
\end{equation*}%
The result now follows, since $u_{k}\leq v_{k_{0}}$ whenever $k\geq k_{0}$.
\end{proof}

\bigskip

\begin{proof}[Proof of Theorem \protect\ref{two}]
Now let $f$, $h$, $\Omega (it_{j})$ and $(S_{m_{k}})$ be as in the statement
of Theorem \ref{two}, and let $0<\varepsilon <(t_{1}-t_{2})/2$. Then $h$ has
a Poisson integral representation of the form 
\begin{equation*}
h(s)=\frac{\sigma }{\pi }\int_{\mathbb{R}}\left\vert s-i\zeta \right\vert
^{-2}d\mu (\zeta )\text{ \ \ \ }(s\in \mathbb{C}_{+}),
\end{equation*}%
where $\mu $ is a positive measure on $\mathbb{R}$. If we define $I=[t_{1}+%
\frac{\varepsilon }{2},t_{2}-\frac{\varepsilon }{2}]$ and 
\begin{equation*}
h_{1}(s)=\frac{\sigma }{\pi }\left( \int_{\mathbb{R}\backslash I}\frac{1}{%
\left\vert s-i\zeta \right\vert ^{2}}d\mu (\zeta )+\frac{\mu (I)}{\left\vert
s-i(t_{1}+\frac{\varepsilon }{2})\right\vert ^{2}}+\frac{\mu (I)}{\left\vert
s-i(t_{1}-\frac{\varepsilon }{2})\right\vert ^{2}}\right) ,
\end{equation*}%
then $e^{-h_{1}}\left\vert f\right\vert $ is also bounded above, by $d\geq 1$
say, on the sets $\Omega (it_{j})$.

There exists $b\geq 1$ such that $\left\vert S_{m_{k}}\right\vert \leq b$ on
the set $K=[i(t_{1}-\varepsilon ),i(t_{2}+\varepsilon )]$, provided $%
\varepsilon $ has been chosen small enough. By Lemma \ref{L}\ there is a
constant $c>0$ such that%
\begin{equation*}
\log \left( \left\vert S_{m_{k}}\right\vert /b\right) \leq c\lambda
_{m_{k}}G_{\mathbb{C}\backslash K}(-2,\cdot )\text{ \ on \ }\{\sigma
>-1\}\backslash K\text{ \ \ \ }(k\in \mathbb{N}).
\end{equation*}%
Hence there is a constant $c_{1}>0$ such that%
\begin{equation*}
\log \left( \left\vert S_{m_{k}}(s)\right\vert /b\right) \leq c_{1}\lambda
_{m_{k}}\sigma \text{ \ on \ }\Omega (it_{j})\text{ \ \ \ }(j=1,2;k\in 
\mathbb{N})
\end{equation*}%
(cf. Lemma 8.5.1 of \cite{AG}). If we define the subharmonic functions%
\begin{equation*}
u_{k}=\frac{1}{\lambda _{m_{k}+1}}\log \left( \frac{\left\vert
S_{m_{k}}-f\right\vert }{2bde^{h_{1}}}\right) \text{ \ \ \ on }\mathbb{C}_{+}%
\text{\ \ \ }(k\in \mathbb{N}),
\end{equation*}%
then 
\begin{eqnarray*}
u_{k}(s) &\leq &\frac{1}{\lambda _{m_{k}+1}}\log ^{+}\left( \frac{\max
\{\left\vert S_{m_{k}}(s)|,|f(s)\right\vert \}}{bde^{h_{1}(s)}}\right) \\
&\leq &\frac{1}{\lambda _{m_{k}+1}}\log ^{+}\max \left\{ \frac{\left\vert
S_{m_{k}}(s)\right\vert }{b},\frac{e^{-h_{1}(s)}\left\vert f(s)\right\vert }{%
d}\right\} \\
&\leq &c_{1}\sigma \text{ \ \ \ }(s\in \Omega (it_{1})\cup \Omega
(it_{2}),k\in \mathbb{N}).
\end{eqnarray*}%
Since we know from Lemma \ref{L0} that $\lim \sup_{k\rightarrow \infty
}u_{k}(s)\leq -\sigma <0$ on $\mathbb{C}_{+}$, it follows from Proposition %
\ref{P} that there exists $k_{0}\in \mathbb{N}$ such that $u_{k}(s)<0$ on
the line segments $(it_{j},it_{j}+1)$ $(j=1,2)$ when $k\geq k_{0}$. Hence $%
\left\vert S_{m_{k}}-f\right\vert \leq 2bde^{h_{1}}$, and so $\log
\left\vert S_{m_{k}}\right\vert \leq \log (2bd+d)+h_{1}$ on these line
segments whenever $k\geq k_{0}$. Since $(S_{m_{k}})$ is uniformly bounded on 
$[it_{1},it_{2}]\cup \lbrack it_{1}+1,it_{2}+1]$, it follows from the
maximum principle that there exists $c_{2}>0$ such that $\log \left\vert
f\right\vert \leq c_{2}+h_{1}$ on the rectangle with vertices $%
\{it_{j},it_{j}+1:j=1,2\}$, and the result follows from the boundedness of $%
h_{1}$ near points of the line segment $[i(t_{1}+\varepsilon
),i(t_{2}-\varepsilon )]$.
\end{proof}

\section{Proof of Theorem \protect\ref{three}}

Let $f(s)=\sum\limits_{n=1}^{\infty }a_{n}e^{-\lambda _{n}s}$, where $f\in D(%
\mathbb{C}_{+})$, and let $\zeta =it_{0}$. We assume that $\left\vert
f^{\prime }\right\vert \leq b$ on some fat approach region $\Omega (it_{0})$
and that $\left\vert S_{m_{k}}\right\vert \leq b$ on the line segment $%
K=[i(t_{0}-\varepsilon ),i(t_{0}+\varepsilon )]$ for some $b\geq 1$ and $%
\varepsilon \in (0,1)$. We know from Lemma \ref{L} that there is a constant $%
c>0$ such that 
\begin{equation*}
\log \left( \left\vert S_{m_{k}}\right\vert /b\right) \leq c\lambda
_{m_{k}}G_{\mathbb{C}\backslash K}(-2,\cdot )\text{ \ on \ }\{\sigma
>-1\}\backslash K\text{ \ \ }(k\in \mathbb{N}).
\end{equation*}%
We choose $k_{0}$ large enough so that $\lambda _{m_{k}}^{-1}<\varepsilon /3$
when $k\geq k_{0}$. Since there is a constant $c_{1}>0$ such that 
\begin{equation*}
G_{\mathbb{C}\backslash K}(-2,s)\leq c_{1}\left\vert \sigma \right\vert 
\text{ \ \ \ }(\left\vert \sigma \right\vert <\varepsilon /3,\left\vert
t-t_{0}\right\vert <2\varepsilon /3),
\end{equation*}%
we see that%
\begin{equation*}
\log \left( \frac{\left\vert S_{m_{k}}(s)\right\vert }{b}\right) \leq cc_{1}%
\text{ \ \ \ }(\left\vert \sigma \right\vert \leq \frac{1}{\lambda _{m_{k}}}%
,\left\vert t-t_{0}\right\vert \leq \frac{\varepsilon }{3}+\frac{1}{\lambda
_{m_{k}}},k\geq k_{0}).
\end{equation*}%
It follows from Cauchy's integral formula for derivatives, applied to
circles of radius $\lambda _{m_{k}}^{-1}$ centred on $L$, that 
\begin{equation*}
\left\vert S_{m_{k}}^{\prime }\right\vert \leq c_{2}\lambda _{m_{k}}\text{ \
\ on \ }L=[i(t_{0}-\varepsilon /3),i(t_{0}+\varepsilon /3)]\text{ \ \
whenever }k\geq k_{0},
\end{equation*}%
where $c_{2}=be^{cc_{1}}$.

Since $f^{\prime }\in D(\mathbb{C}_{+})$ (see Theorem 4 of \cite{HR}), we
can appeal to Lemma \ref{L} once again to see that 
\begin{equation*}
\log \left( \frac{\left\vert S_{m_{k}}^{\prime }\right\vert }{c_{2}\lambda
_{m_{k}}}\right) \leq c_{3}\lambda _{m_{k}}G_{\mathbb{C}\backslash
L}(-2,\cdot )\text{ \ on }\mathbb{C}_{+}\text{ \ \ \ }(k\geq k_{0}),
\end{equation*}%
for some constant $c_{3}>0$. We define the subharmonic functions 
\begin{equation*}
u_{k}=\frac{1}{\lambda _{m_{k}+1}}\log \left( \frac{\left\vert
S_{m_{k}}^{\prime }-f^{\prime }\right\vert }{2c_{2}\lambda _{m_{k}}}\right) 
\text{ \ on \ }\mathbb{C}_{+}\text{ \ \ \ }(k\in \mathbb{N}),
\end{equation*}%
and note from Lemma \ref{L0} that $\lim \sup_{k\rightarrow \infty
}u_{k}(s)\leq -\sigma <0$ on\ $\mathbb{C}_{+}$. Further,%
\begin{equation*}
u_{k}\leq \frac{1}{\lambda _{m_{k}+1}}\log \left( \frac{\max \{\left\vert
S_{m_{k}}^{\prime }|,|f^{\prime }\right\vert \}}{c_{2}\lambda _{m_{k}}}%
\right) \leq c_{3}G_{\mathbb{C}\backslash L}(-2,\cdot )\text{ \ on \ }\Omega
(it_{0})\text{ \ \ }(k\geq k_{0}).
\end{equation*}%
Since $G_{\mathbb{C}\backslash L}(-2,s)$ is comparable to $\left\vert \sigma
\right\vert $ near $it_{0}$, we see from Proposition \ref{P} that there
exist $k_{1}\in \mathbb{N}$ and $c_{4}>0$ such that 
\begin{equation*}
\left\vert S_{m_{k}}^{\prime }-f^{\prime }\right\vert (s)\leq 2c_{2}\lambda
_{m_{k}}\exp (-c_{4}\lambda _{m_{k}+1}\sigma )\text{ \ \ \ }(s\in
(it_{0},\xi _{0}],k\geq k_{1}),
\end{equation*}%
where $\xi _{0}=it_{0}+1$. Hence%
\begin{eqnarray*}
\left\vert \;\left\vert S_{m_{k}}(it_{0})-f(it_{0})\right\vert -\left\vert
S_{m_{k}}(\xi _{0})-f(\xi _{0})\right\vert \;\right\vert &\leq
&\int_{[it_{0},\xi _{0}]}\left\vert S_{m_{k}}^{\prime }(s)-f^{\prime
}(s)\right\vert \left\vert ds\right\vert \\
&\leq &2c_{2}\lambda _{m_{k}}\int_{0}^{1}\exp (-c_{4}\lambda
_{m_{k}+1}\sigma )d\sigma \\
&\leq &2\frac{c_{2}}{c_{4}}\frac{\lambda _{m_{k}}}{\lambda _{m_{k}+1}}%
\rightarrow 0\text{ \ \ \ }(k\rightarrow \infty ),
\end{eqnarray*}%
by hypothesis. Since $S_{m_{k}}(\xi _{0})\rightarrow f(\xi _{0})$ we also
have $S_{m_{k}}(it_{0})\rightarrow f(it_{0})$, as required.

\section{Proofs of Applications\label{pfapps}}

\begin{proof}[Proof of Corollary \protect\ref{C1}]
Let $f$ be a holomorphic function on $\mathbb{C}_{+}$. Plessner's theorem
tells us that, at almost every point $\zeta $ of $i\mathbb{R}$, either $f$
has a finite nontangential limit or $f(\Gamma )$ is dense in $\mathbb{C}$
for every open isosceles triangle $\Gamma \subset \mathbb{C}_{+}$ with
vertex at $\zeta $ and horizontal axis of symmetry. Theorem \ref{one} tells
us that $f$ cannot both belong to $UD_{0}$ and also have a finite
nontangential limit on a subset of $i\mathbb{R}$ of positive measure. Thus
Corollary \ref{C1} is proved.
\end{proof}

\bigskip

\begin{proof}[Proof of Corollary \protect\ref{C2}]
Let $f\in UD_{0}$. To see why Corollary \ref{C2} holds, suppose (for the
sake of contradiction) that there is a disc $D$, centred at a point of $i%
\mathbb{R}$, for which the set $\mathbb{C}\backslash f(D\cap \mathbb{C}_{+})$
has positive logarithmic capacity. Myrberg's theorem (Theorem 5.3.8 of \cite%
{AG}) then shows that there is a positive harmonic function $h$ on the open
set $f(D\cap \mathbb{C}_{+})$ satisfying $h(z)\geq \log \left\vert
z\right\vert $ there. Thus $\log \left\vert f\right\vert \leq h\circ f$ on $%
D\cap \mathbb{C}_{+}$. Fatou's theorem now tells us that $h\circ f$, and
hence $f$, is nontangentially bounded at almost every point of $D\cap i%
\mathbb{R}$, contradicting Corollary \ref{C1}.
\end{proof}

\bigskip

\begin{proof}[Proof of Corollary \protect\ref{C3}]
Let $f\in UD_{0}$. By Corollary \ref{C2}, the unrestricted cluster set of $f$
at each boundary point of $\mathbb{C}_{+}$ is all of $\mathbb{C}$. The
Collingwood Maximality Theorem (Theorem 4.8 of \cite{CL}), when adapted to
the half-plane, tells us that there is a residual set $Z\subset \mathbb{R}$
such that, for each $t\in Z$, the cluster set of $f$ along $\{\sigma
+it:0<\sigma <1\}$ as $\sigma \rightarrow 0$ is also all of $\mathbb{C}$.
The result now follows.
\end{proof}

\bigskip

\begin{proof}[Proof of Corollary \protect\ref{C4}]
Let $f\in D(\mathbb{C}_{+})$ and $h$ be a positive harmonic function on $%
\mathbb{C}_{+}$, and suppose that $\left\vert f\right\vert \leq e^{h}$ on
two fat approach regions $\Omega (it_{j})$, where $j=1,2$ and $t_{1}<t_{2}$.
If $f\in UD_{0}$, then we can find a subsequence $(S_{m_{k}})$ of partial
sums converging uniformly to $0$ on $[i(t_{1}-1),i(t_{2}+1)]$. Theorem \ref%
{two} and Remark \ref{R0} then show that $f$ is bounded near all points of
the interval $(it_{1},it_{2})$ and has nontangential limit $0$ almost
everywhere on this line segment. This leads to the contradictory conclusion
that $f\equiv 0$.
\end{proof}

\bigskip

In what follows the partial sums of a power series $\sum_{0}^{\infty
}a_{n}z^{n}$ will be denoted by $(T_{m})$. The series is said to have pure
Ostrowski gaps $(p_{k},q_{k})$ if 
\begin{equation*}
1\leq p_{1}<q_{1}\leq p_{2}<q_{2}\leq ...\text{ \ and }q_{k}/p_{k}%
\rightarrow \infty ,
\end{equation*}%
and if 
\begin{equation*}
a_{n}=0\text{ \ whenever \ }p_{k}+1\leq n\leq q_{k}\text{ \ for some \ }k\in 
\mathbb{N}.
\end{equation*}

The following consequence of Theorem \ref{three} is an internediate step
towards proving Corollary \ref{C5}.

\begin{corollary}
\label{C6}Let $f$ be a holomorphic function on $\mathbb{D}$ such that $%
f^{\prime }$ is bounded on a fat approach region to some point $w$ of $%
\partial \mathbb{D}$ (whence $f(w)$ exists as a nontangential limit). If the
Taylor series of $f$ about $0$ has pure Ostrowski gaps $(p_{k},q_{k})$, and $%
(T_{p_{k}})$ is uniformly bounded on an open arc of $\partial \mathbb{D}$
containing $w$, then $T_{p_{k}}(w)\rightarrow f(w)$.
\end{corollary}

\bigskip

\begin{proof}[Proof of Corollary \protect\ref{C6}]
Without loss of generality we may assume that $w=1$. Then the function $%
g:s\longmapsto f(e^{-s})$ belongs to $D(\mathbb{C}_{+})$ and has bounded
derivative on a fat approach region to $0\in i\mathbb{R}$. Since $g$ has the
form%
\begin{equation*}
g(s)=\sum_{k=1}^{\infty }\sum_{n=q_{k-1}+1}^{p_{k}}a_{n}e^{-ns}\text{ \ \ \ }%
(s\in \mathbb{C}_{+}),
\end{equation*}%
where $q_{0}$ is interpreted as $-1$, and $(q_{k}+1)/p_{k}\rightarrow \infty 
$, the result follows from Theorem \ref{three}.
\end{proof}

\bigskip

\begin{proof}[Proof of Corollary \protect\ref{C5}]
Let $f(z)=\sum_{0}^{\infty }a_{n}z^{n}$ be a holomorphic function on $%
\mathbb{D}$ such that $f^{\prime }$ is bounded on a fat approach region to
some point $e^{i\phi }$ of $\partial \mathbb{D}$. Now suppose, for the sake
of contradiction, that $f\in U(\mathbb{D},0)$. It follows from Theorem 2 of 
\cite{MVY} that we can find sequences $(p_{k})$ and $(q_{k})$ such that

\begin{enumerate}
\item[(i)] $f=g+h$ on $\mathbb{D}$, where $h$ is entire and $g$ possesses
pure Ostrowski gaps $(p_{k},q_{k})$; and

\item[(ii)] $T_{p_{k}}\rightarrow f(e^{i\phi })+1$ \ uniformly on $%
\{e^{i\theta }:\left\vert \theta -\phi \right\vert \leq \pi /2\}$.
\end{enumerate}

We can now apply Corollary \ref{C6} to the function $f-h$ to see that $%
T_{p_{k}}(e^{i\phi })\rightarrow f(e^{i\phi })$, which contradicts (ii).
\end{proof}

\bigskip

\noindent School of Mathematics and Statistics,

\noindent University College Dublin,

\noindent Belfield, Dublin 4, Ireland.

\noindent e-mail: stephen.gardiner@ucd.ie

\bigskip

\noindent Department of Mathematics,

\noindent University of Western Ontario,

\noindent London, Ontario, Canada N6A 5B7.

\bigskip

\noindent \textit{Current address for Myrto Manolaki:}

\noindent School of Mathematics and Statistics,

\noindent University College Dublin,

\noindent Belfield, Dublin 4, Ireland.

\noindent e-mail: arhimidis8@yahoo.gr

\end{document}